\documentclass{amsart}

\usepackage{amssymb,graphicx}
\usepackage{enumerate}

\usepackage[all]{xy}
\usepackage{epstopdf}

\newtheorem{theorem}{Theorem}
\newcommand{\bt}{\begin{theorem}}
\newcommand{\et}{\end{theorem}}
\newtheorem{lemma}{Lemma}
\newcommand{\bl}{\begin{lemma}}
\newcommand{\el}{\end{lemma}}
\newtheorem{corollary}{Corollary}
\newcommand{\bc}{\begin{corollary}}
\newcommand{\ec}{\end{corollary}}
\newtheorem{problem}{Problem}
\newcommand{\bprob}{\begin{problem}}
\newcommand{\eprob}{\end{problem}}
\newcommand{\beq}{\begin{equation}}
\newcommand{\eeq}{\end{equation}}
\newcommand{\benum}{\begin{enumerate}}
\newcommand{\eenum}{\end{enumerate}}

\newcommand{\Z}{\ensuremath{\mathbf Z}}

\newcommand{\R}{\ensuremath{\mathbf R}}

\newcommand{\mci}{\ensuremath{ \mathcal I}}

\DeclareMathOperator{\id}{id}

\newcommand{\bmat}{\left(\begin{matrix}}
\newcommand{\emat}{\end{matrix}\right)}

\DeclareMathOperator{\qqand}{\qquad\text{and}\qquad}

\title{MSTD sets and Freiman isomorphisms}
\author{Melvyn B. Nathanson}
\address{Department of Mathematics\\
Lehman College (CUNY)\\
Bronx, NY 10468} 
\email{melvyn.nathanson@lehman.cuny.edu}

\dedicatory{To Grigori Freiman on his 90th birthday}


\date{\today}

\begin{document}
\maketitle

\begin{abstract}
An MSTD set is a finite set with more pairwise sums than differences.  
$(\Upsilon,\Phi)$-ismorphisms are generalizations of Freiman isomorphisms  
to arbitrary linear forms.  These generalized isomorphisms are used to 
prove that every finite set of real numbers is Freiman isomorphic to a finite set 
of integers.  This implies that there exists no MSTD set $A$ of real numbers 
with $|A| \leq 7$, and, up to Freiman isomorphism and affine isomorphism, there exists 
exactly one MSTD set $A$ of real numbers with $|A| = 8$.
\end{abstract}

\section{Sums and differences}

For every nonempty subset $A$ of an additive abelian group, we define the \textit{sumset} 
\[
A+A = \{a+a': a,a' \in A\}
\]
and the \textit{difference set} 
\[
A - A = \{a- a': a,a' \in A\}.
\]
The observation that 
\[
3+2 = 2 + 3 
\]
but
\[
3-2 \neq 2-3
\]
suggests the reasonable conjecture that a finite set of integers or real numbers 
should have more differences, or at least as many differences, as sums, 
but this conjecture is false.  For example, the set 
\[
A = \{0, 2, 3, 4, 7, 11, 12, 14\}
\]
has  difference set 
\[
A - A = [-14,14] \setminus \{  6,-6, 13, -13 \}
\]
and sumset  
\[
A+A = [0,28] \setminus \{1,20,27\}.  
\]
Note the translated difference set 
\[
A-A + \{ 14\} =  [0,28] \setminus \{1, 8, 20,27\}.  
\]
We have  
\[
 |A -A| = 25 < 26 = |A+A|.
\]
Sets with more sums than differences are called  \textit{MSTD sets}.  

As expected, MSTD sets of integers are rare (e.g. Hegarty and Miller~\cite{hega-mill09},  
Martin and O'Bryant~\cite{mart-obry07a}, Zhao~\cite{zhao10b,zhao11}).  
Isolated examples and infinite families of MSTD sets of integers have been constructed 
(e.g.  Hegarty~\cite{hega07a}, 
Miller, Orosz, and Scheinerman~\cite{mill-oros-sche10}, 
Miller and Scheinerman~\cite{mill-sche10}, 
and Nathanson~\cite{nath07xyz,nath07yzx}), 
but there is no adequate classification.  

If $A$ is an MSTD set, then every affine image of $A$, that is, every set of the form 
\[
\lambda \ast A + \mu = \{ \lambda a + \mu: a \in A\}
\]
with $\lambda, \mu \in \R$ and $\lambda \neq 0$ is also an MSTD set.  
We call the sets $A$ and $\lambda \ast A + \mu $  \textit{affinely isomorphic}.  
Using a clever algorithm and extensive computation, 
Hegarty~\cite{hega07a} proved that if $A$ is an MSTD set of integers, 
then $|A| \geq 8$, and if $A$ is an MSTD set with $|A| = 8$, 
then $A$ is affinely isomorphic to the set 
$\{0, 2, 3, 4, 7, 11, 12, 14\}$.

MSTD sets of integers have been studied extensively, but little is known 
about MSTD sets of real numbers and how they differ from MSTD sets of integers.  
In this paper we show that, associated to every MSTD set of real numbers, 
there is a ``Freiman isomorphic'' MSTD set of integers,  
and that properties of MSTD sets of real numbers may be 
deduced from properties of their associated  integer MSTD sets.  
In particular, we prove that if $A$ is an MSTD set of real numbers, 
then $|A| \geq 8$.  We shall define \textit{Freiman isomorphism} 
and, more generally, \textit{$(\Upsilon,\Phi)$-isomorphism}, 
and prove that if $A$ is a real MSTD set with $|A| = 8$, then $A$ is 
not only Freiman isomorphic but also affinely isomorphic to 
the set $\{0, 2, 3, 4, 7, 11, 12, 14\}$.

\section{$( \Upsilon, \Phi)$-isomorphisms} 

We begin with a simple set theoretical observation.  
Let $U$ be a set, and let $\Upsilon$ be a function defined on $U$ 
with image $\Upsilon(U) = \{ \Upsilon(u): u \in U\}$.
Let $V$ be a set, and let $\Phi$ be a function defined on $V$ 
with image $\Phi(V) = \{ \Phi(v): v \in V \}$.
A function $f^*: U \rightarrow V$ is an \textit{$(\Upsilon, \Phi)$-homomorphism} 
if, for all $u, u' \in U$, the equation   
\beq          \label{MPTQ:Upsilon-f}
\Upsilon(u) =  \Upsilon(u') 
\eeq
implies that  
\beq          \label{MPTQ:Phi-f}
\Phi(f^* (u)) =  \Phi(f^* (u')).
\eeq
Equivalently, for each $x \in \Upsilon(U)$, there exists $y \in \Phi(V)$ 
such that 
\[
f^*\left(  \Upsilon^{-1}(x) \right) \subseteq \Phi^{-1}(y).
\]
This implies that there is a well-defined function $F: \Upsilon(U)  \rightarrow \Phi(V)$ 
such that 
\[
F(x) = \Phi(f^* (u))
\]
for all $x \in \Upsilon(U)$ and $u \in \Upsilon^{-1}(x)$.  
Thus, the following diagram commutes:
\[
\xymatrix{
U \ar[r]^{\Upsilon} \ar[d]_{f^*}  &\Upsilon(U)  \ar[d]^F  \\
V \ar[r]^{\Phi}  &  \Phi(V)
}
\]

An $(\Upsilon, \Phi)$-homomorphism  $f^*: U\rightarrow V$ is a 
\textit{$(\Upsilon, \Phi)$-isomorphism} 
if there is a  $(\Phi, \Upsilon)$-homomorphism $g^*: V \rightarrow U$ such that 
$g^* f^* = \id_U$ and $f^* g^* = \id_V$.  
We call $g^*$ the \textit{$(\Upsilon, \Phi)$-inverse} of $f^*$.  
As a $(\Phi, \Upsilon)$-homomorphism, 
the function $g^*$ has the following property:  
If $v, v' \in V$ and   
\beq     							     \label{MPTQ:Phi-g}
\Phi(v) = \Phi(v')
\eeq
then  
\beq       							   \label{MPTQ:Upsilon-g}
\Upsilon( g^*(v)) =  \Upsilon( g^*(v')).
\eeq

Let $f^*$ be an $(\Upsilon, \Phi)$-isomorphism with inverse $g^*$.  
Let $G: \Phi(V) \rightarrow \Upsilon(U)$ be the function induced by $g^*$.  
Let $u, u' \in U$ satisfy~\eqref{MPTQ:Phi-f}.  
Writing $v = f^*(u)$ and $v' = f^*(u')$, we obtain~\eqref{MPTQ:Phi-g}, 
which implies~\eqref{MPTQ:Upsilon-g}.
Because $u = g^*(v)$ and $u' = g^*(v')$, we obtain~\eqref{MPTQ:Upsilon-f}. 
Thus, if $f^*$ is an $(\Upsilon, \Phi)$-isomorphism, then~\eqref{MPTQ:Upsilon-f}
and~\eqref{MPTQ:Phi-f} are equivalent.

\bl
Let $U, V,$ and $W$ be sets, and let $\Upsilon, \Phi,$ and $\chi$ 
be functions on $U, V,$ and $W$, respectively.  
If $f^*$ is an $(\Upsilon, \Phi)$-homomorphism  
(resp. $(\Upsilon, \Phi)$-isomorphism) from $U$ to $V$ 
and if $f_1^*$ is a $(\Phi, \chi)$-homomorphism 
(resp. $(\Phi, \chi)$-isomorphism) from $V$ to $W$, then 
$f_1^* f^*$ is an $(\Upsilon, \chi)$-homomorphism
(resp. $(\Upsilon, \chi)$-isomorphism) from $U$ to $W$.
\el

\begin{proof}
This  follows immediately from the definitions.  
\end{proof}

Let $|X|$ denote the cardinality of the set $X$.  
For every set $U$ and function $\Upsilon$ on $U$, 
we define the \textit{representation function} of $x \in \Upsilon(U)$ as follows:
\[
r_{\Upsilon,U}(x) 
= \left| \left\{u \in U : \Upsilon (u)  = x \right\} \right|
=  \left| \Upsilon^{-1} (x) \right|.  
\]

\bt                  \label{MPTQ:theorem:UVrep}
Let $U$ and $V$ be sets, let $\Upsilon$ be a function on $U$, 
and let $\Phi$ be a function on $V$. 
If  $f^*:U \rightarrow V$ is an $(\Upsilon, \Phi)$-isomorphism, then 
\benum
\item[(i)]
The induced function $F: \Upsilon(U) \rightarrow \Phi(V)$ is a bijection, 
and so 
\[
|\Upsilon(U)| = |\Phi(V)|.
\]
\item[(ii)]
For all $x \in \Upsilon(U)$, 
\[
r_{\Upsilon,U}(x) = r_{\Phi,V}(F(x)). 
\]
\eenum
\et

\begin{proof}
Let  $g^*$ be the $(\Upsilon, \Phi)$-inverse of $f^*$.  
If $x \in \Upsilon(U)$ and $u \in \Upsilon^{-1}(x)$, then 
\[
GF(x)  = G F \Upsilon(u) = G \Phi f^* (u)  = \Upsilon g^* f^* (u) = \Upsilon(u) = x.
\]
Similarly, $FG(y) = y$ for all $y \in \Phi(V)$, and so $G = F^{-1}$.  
Thus, $F$ is a bijection and $|\Upsilon(U)| = |\Phi(V)|$.  

For all $x \in \Upsilon(U)$, we have 
\[
u \in \Upsilon^{-1}(x) 
\]
if and only if 
\[
F\Upsilon(u) = F(x)
\]
if and only if 
\[
\Phi f^*(u) = F(x) 
\]
if and only if 
\[
f^*(u) \in \Phi^{-1}( F(x) ).
\]
Because $f^*:U \rightarrow V$ is a bijection, 
we have $|\Upsilon^{-1}(x)| = |\Phi^{-1}( F(x) ) |$, or, equivalently, 
$r_{\Upsilon,U}(x) = r_{\Phi,V}(F(x))$.  
This completes the proof.  
\end{proof}

\section{Linear forms}
A \textit{linear form} with coefficients in a ring $R$ is a function of the form 
\[
\Upsilon(t_1,\ldots, t_h) = \sum_{j=1}^h \upsilon_j t_j
\]
where $\upsilon_j \in R$ for $j = 1, \ldots, h$. 
The linear form is nonzero if $\upsilon_j \neq 0$ for some $j$.  
The study of sums and differences is the study of the special linear forms 
$t_1+t_2$ and $t_1 - t_2$, but it is natural to consider more general linear forms.  
We begin with functions of $h$ variables that are sums of functions of one variable.

\bt                                    \label{MPTQ:theorem:UpsilonPhi-J}
Let $R$ be a ring, let $M$ be an $R$-module, let $A$ be a subset of  $M$, 
and let $A^h$ denote the set of all $h$-tuples of elements of $A$.
Let $S$ be a ring, let $N$ be an $S$-module,  let $B$ be a subset of  $N$, 
and let $B^h$ denote the set of all $h$-tuples of elements of $B$.  
For every function $f:A \rightarrow B$, define 
$f^*:A^h \rightarrow B^h$ by 
\[
f^*(a_1,\ldots, a_h) = (f(a_1), \ldots, f(a_h)).
\] 
Let $\Upsilon_1(t_1), \ldots, \Upsilon_h(t_h)$ be functions defined on $A$, 
and  let $\Phi_1(t_1), \ldots, \Phi_h(t_h)$ be functions defined on $B$.  
The functions 
\[
\Upsilon(t_1,\ldots, t_h) = \sum_{j=1}^h \Upsilon_j (t_j)
\]
and  
\[
\Phi(t_1,\ldots, t_h) = \sum_{j=1}^h \Phi_j (t_j) 
\]
are defined on $A^h$ and $B^h$, respectively.

For every subset $J$ of $H = \{1,2,\ldots, h\}$, let 
\[
\Upsilon_J(t_1,\ldots, t_h) = \sum_{j \in H\setminus J} \Upsilon_j (t_j) 
-  \sum_{j \in J} \Upsilon_j (t_j)  
\]
and
\[
\Phi_J(t_1,\ldots, t_h) = \sum_{j \in H\setminus J} \Phi_j (t_j)
-  \sum_{j \in J}  \Phi_j (t_j).   
\]
The functions $\Upsilon_J$ and $\Phi_J$ are functions on $A^h$ 
and $B^h$, respectively.  
If $f^*:A^h \rightarrow B^h$ is an $(\Upsilon,\Phi)$-isomorphism,  
then, for every subset $J \subseteq H$, the function $f^*:A^h \rightarrow B^h$ 
is also an $(\Upsilon_J,\Phi_j)$-isomorphism, and 
\[
\left| \Upsilon_J(A^h) \right| = \left| \Phi_J(B^h) \right|.  
\] 
\et

\begin{proof}
If $(a_1,\ldots, a_h) \in A^h$ and $(a'_1,\ldots, a'_h) \in A^h$ 
and if 
\[
\Upsilon_J (a_1,\ldots, a_h) = \Upsilon_J (a'_1,\ldots, a'_h) 
\]
then 
\[
\sum_{j \in H\setminus J} \Upsilon_j (a_j) -  \sum_{j \in J} \Upsilon_j (a_j)  
= \sum_{j \in H\setminus J} \Upsilon_j (a'_j) -  \sum_{j \in J} \Upsilon_j (a'_j)  
\]
and so 
\[
\sum_{j \in H\setminus J} \Upsilon_j (a_j) +  \sum_{j \in J} \Upsilon_j (a'_j)  
= \sum_{j \in H\setminus J} \Upsilon_j (a'_j)  + \sum_{j \in J} \Upsilon_j (a_j). 
\]
Let
\[
u_j = \begin{cases}
a_j & \text{if $j \in H \setminus J$} \\
a'_j & \text{if $j \in J$} 
\end{cases}
\]
and
\[
u'_j = \begin{cases}
a'_j & \text{if $j \in H \setminus J$} \\
a_j & \text{if $j \in J$.} 
\end{cases}
\]
We have $(u_1,\ldots, u_h), (u'_1,\ldots, u'_h) \in A^h$, and 
\[
\Upsilon(u_1,\ldots, u_h)  = \Upsilon(u'_1,\ldots, u'_h). 
\]
Because $f^*$ is an $(\Upsilon,\Phi)$-homomorphism, it follows that 
\[
\Phi(   f^*( u_1,\ldots, u_h)   ) = \Phi(  f^*( u'_1,\ldots, u'_h) ). 
\]
Expanding the left and right sides of this equation, we obtain 
\begin{align*}
\Phi(   f^*( u_1,\ldots, u_h)   ) 
& = \Phi( f( u_1),\ldots,  f( u_h) ) \\
& =  \sum_{j \in H \setminus J} \Phi_j (f( u_j)) +  \sum_{j \in J}  \Phi_j (f( u_j)) \\  
& =  \sum_{j \in H \setminus J} \Phi_j (f( a_j)) +  \sum_{j \in J}  \Phi_j (f( a'_j)) 
\end{align*}
and
\begin{align*}
\Phi( f^*( u'_1,\ldots, u'_h) ) 
& = \Phi( f( u'_1),\ldots,  f( u'_h) ) \\
& =  \sum_{j \in H \setminus J} \Phi_j (f( u'_j)) + \sum_{j \in J}  \Phi_j (f( u'_j)) \\  
& =  \sum_{j \in H \setminus J} \Phi_j (f( a'_j)) + \sum_{j \in J}  \Phi_j (f( a_j)).  
\end{align*}
Therefore,
\begin{align*}
 \Phi_J(  f^*( a_1, \ldots, a_h) )  
&  =  \sum_{j \in H \setminus J} \Phi_j (f( a_j)) -  \sum_{j \in J}  \Phi_j (f( a_j)) \\
& =  \sum_{j \in H \setminus J} \Phi_j (f( a'_j))  -  \sum_{j \in J}  \Phi_j (f( a'_j)) \\
& = \Phi_J( f^*( a'_1, \ldots, a'_h) ).
 \end{align*}
Similarly, because $\left(f^*\right)^{-1}$ is a $(\Phi,\Upsilon)$-homomorphism,
\[
\Phi_J (f^*( a_1, \ldots, a_h) ) = \Phi_J (f^*( a'_1, \ldots, a'_h) )
\]
implies that 
\[
\Upsilon_J (a_1,\ldots, a_h) = \Upsilon_J (a'_1,\ldots, a'_h) 
\]
and so $f^*$ is also an $(\Upsilon_J,\Phi_J)$-isomorphism for all $J \subseteq H$.  
Applying Theorem~\ref{MPTQ:theorem:UVrep} completes the proof.  
\end{proof}

\section{Linear forms with rational coefficients}
In this section we consider the special case $\Upsilon = \Phi$, 
and write, simply,  \textit{$\Phi$-homomorphism}  
instead of  $(\Phi, \Phi)$-homomorphism,  
and  \textit{$\Phi$-isomorphism} instead of $(\Phi, \Phi)$-isomorphism.

Let $A = \{a_1,\ldots, a_k\}$ and $B = \{b_1,\ldots, b_k\}$ 
be finite sets of real numbers with $|A| = |B| = k \geq 2$, and let $f:A \rightarrow B$ 
be a bijection such that $f(a_i) = b_i$ for all $i = 1, \ldots, k$.  
Let $\Phi(t_1,\ldots, t_h)$ be a nonzero linear form with rational coefficients.  
The function $f$ will be called a $\Phi$-isomorphism if the function 
$f^*:A^h \rightarrow B^h$ is a $\Phi$-isomorphism.  
This means that, for all $h$-tuples 
$(a_{i_1}, \ldots, a_{i_h}), (a_{i_{h+1}} , \ldots, a_{i_{2h}} )  \in A^h$, 
we have 
\[
\Phi(a_{i_1}, \ldots, a_{i_h} ) = \Phi( a_{i_{h+1}} , \ldots, a_{i_{2h}} ) 
\]
if and only if 
\[
\Phi(b_{i_1}, \ldots, b_{i_h} ) = \Phi( b_{i_{h+1}} , \ldots, b_{i_{2h}}).  
\]
The sets $A$ and $B$ are \textit{$\Phi$-isomorphic} if there exists 
a $\Phi$-isomorphism $f:A \rightarrow B$.  
The following result generalizes Corollary 8.1 of Nathanson~\cite{nath96bb}.

\bt               \label{MPTQ:theorem:lattice}
Let 
\[
\Phi(t_1,\ldots, t_h) = \sum_{j =1}^h \varphi_j  t_j
\]
be a nonzero linear form with rational coefficients.  
Every finite subset of $\Z^d$ is $\Phi$-isomorphic to a set of positive integers.  
\et

\begin{proof}
Let $m$ be a common multiple of the denominators of the rational numbers 
$\varphi_1, \ldots, \varphi_h$.  The linear form
\[
\hat{\Phi} (t_1,\ldots, t_h) = m\Phi(t_1,\ldots, t_h) = \sum_{j =1}^h m\varphi_j t_j
\]
has integer coefficients.  For all $h$-tuples of lattice points or real numbers, 
we have 
$\hat{\Phi} (a_{1}, \ldots, a_{h} ) = \hat{\Phi} ( a_{h+1} , \ldots, a_{2h}) $ 
if and only if 
$\Phi(a_{1}, \ldots, a_{h} ) = \Phi( a_{h+1} , \ldots, a_{2h} )$.
Thus, we can assume that the coefficients $\varphi_1, \ldots, \varphi_h$ are integers.

For $i =1,\ldots, d$, let $\pi_i: \R^d \rightarrow \R$ 
be the projection of a vector onto its $i$th coordinate.  
Let $A$ be a nonempty finite subset of $\Z^d$, and let 
\[
a^* = \max(\{ |\pi_i(a)| : a \in A \text{ and } i=1,\ldots, d ) \}.
\]
Let 
\[
\varphi^* = \max \left( \{|\varphi_j| : j =1,\ldots, h \} \right) 
\]
and let $\lambda$ be a real number such that 
\beq   \label{MPTQ:ineqT}
\lambda > 2a^* \varphi^* h + 1
\eeq
We shall prove that the linear form $f_{\lambda}: \R^d \rightarrow \R$  defined by 
\[
f_{\lambda}(x_1,\ldots, x_d) = \sum_{i =1}^d  x_i \lambda^{i - 1}
\]
is a $\Phi$-isomorphism from $A$ into \R.  
If $\lambda$ is an integer, then $f$ is a $\Phi$-isomorphism from $A$ into \Z.

Let $(a_1, \ldots, a_h)$ and $ (a_{h+1}, \ldots, a_{2h})$ be $h$-tuples of lattice points in $A$.  
The identity  
\begin{align*}
\Phi(a_1,\ldots, & a_h)  - \Phi(a_{h+1},\ldots, a_{2h}) 
 = \sum_{j=1}^h \varphi_j ( a_j - a_{h+j})    \\
& = \left(  \sum_{j=1}^h \varphi_j \pi_1( a_j - a_{h+j}), \ldots,  \sum_{j=1}^h \varphi_j \pi_d( a_j - a_{h+j}  )  \right)  \in \Z^d
\end{align*}
implies that  
\[
\Phi(a_1,\ldots,  a_h)  = \Phi(a_{h+1},\ldots, a_{2h}) 
\]
if and only if 
\beq   \label{MPTQ:conditionT}
\sum_{j=1}^h \varphi_j \pi_i ( a_j - a_{h+j}) = 0  \qquad\text{ for  all $i = 1,\ldots, d$}.
\eeq
We have 
\[
f_{\lambda}\left( \Phi(a_1,\ldots, a_h) \right) = f_{\lambda}\left(  \Phi(a_{h+1},\ldots, a_{2h})  \right) 
\]
if and only if 
\[
f_{\lambda}\left( \Phi(a_1,\ldots, a_h)  -  \Phi(a_{h+1},\ldots, a_{2h})  \right) = 0
\]
if and only if 
\[
f_{\lambda}\left(  \sum_{j=1}^h \varphi_j \pi_1( a_j - a_{h+j}), \ldots,  \sum_{j=1}^h \varphi_j \pi_d( a_j - a_{h+j}  )  \right) = 0
\]
if and only if 
\beq   \label{MPTQ:eqT}
\sum_{i=1}^d  \left(  \sum_{j =1}^h \varphi_j \pi_i( a_j - a_{h+j})  \right)  \lambda^{i -1} = 0.  
\eeq
Suppose that 
\[
 \sum_{j =1}^h \varphi_j \pi_i( a_j - a_{h+j}) \neq 0
\]
for some $i \in \{1,\ldots, d\}$.  Let $r$ be the greatest integer such that 
$ \sum_{j =1}^h \varphi_j \pi_r( a_j - a_{h+j }) \neq 0$, or, equivalently, such that 
\beq   \label{MPTQ:ineqTT}
 \left|  \sum_{j =1}^h \varphi_j \pi_r( a_j - a_{h+j}) \right| \geq 1.
\eeq
Rewriting~\eqref{MPTQ:eqT}, we obtain 
\[
- \sum_{j =1}^h \varphi_j \pi_r( a_j - a_{h+ j}) \lambda^{r-1} 
= \sum_{i=1}^{r-1}  \left(  \sum_{j =1}^h \varphi_j \pi_i( a_j - a_{h+j})  \right)  \lambda^{i -1}.   
\]
Applying  inequality~\eqref{MPTQ:ineqTT}, the triangle inequality,  
and inequality~\eqref{MPTQ:ineqT}, we obtain 
\begin{align*}
 \lambda^{r-1}  & \leq \left|  \sum_{j =1}^h \varphi_j \pi_r( a_j - a_{h+ j}) \lambda^{r-1}  \right| \\
& = \left|  \sum_{i=1}^{r-1}  \left(  \sum_{j =1}^h \varphi_j \pi_i( a_j - a_{h+j})  \right)  \lambda^{i -1} \right| \\
& \leq  \sum_{ i =1}^{r-1}  \sum_{ j =1}^h \varphi_j  (| \pi_i ( a_j )| + | \pi_i ( a_{h+j} ) |)  \lambda^{i-1} \\
& \leq  \sum_{ i =1}^{r-1}  2 a^*  \varphi^* h \lambda^{ i -1} \\ 
& <  2a^* \varphi^* h  \frac{\lambda^{r-1}}{\lambda - 1} \\
& <  \lambda^{r-1}
\end{align*}
which is absurd.  
Therefore,
\[
f_{\lambda}\left( \Phi(a_1,\ldots, a_h) \right) = f_{\lambda}\left(  \Phi(a_{h+1},\ldots, a_{2h})  \right) 
\]
if and only if condition~\eqref{MPTQ:conditionT} is satisfied. 
that is,  if and only if 
\[
\Phi(a_1,\ldots, a_h) =  \Phi(a_{h+1},\ldots, a_{2h}). 
\]
Thus, $f_{\lambda}$ is a $\Phi$-isomorphism from $A$ to a set of real numbers,   
which, by translation, is  $\Phi$-isomorphic to a set of positive real numbers.  
Choosing a positive integer $\lambda$ completes the proof.   
 \end{proof}

We shall give three different proofs of the following fundamental theorem.
The first uses linear programming, the second uses diophantine approximation, 
and the third reduces the proof to Theorem~\ref{MPTQ:theorem:lattice}.

\bt            \label{MPTQ:theorem:main}
Let
\[
\Phi(t_1,\ldots, t_h) = \sum_{j=1}^h \varphi_j t_j 
\]
be a nonzero linear form with rational coefficients.  
Every nonempty finite set $A$ of real numbers is $\Phi$-isomorphic to a set of positive integers.  
\et

\begin{proof}
There is nothing to prove if  $|A| = 1$, so we can assume that 
$A = \{a_1,\ldots, a_k\}$ is a  finite set of real numbers with $|A| = k \geq 2$.

The first proof uses a linear programming argument of 
Alon and Kleitman~\cite[Prop. 4.1']{alon-klei90}.  
To every pair of $h$-tuples $(a_{i_1}, \ldots, a_{i_h}), (a_{i_{h+1}} , \ldots, a_{i_{2h}} )  \in A^h$, 
we construct a linear equation or linear inequality
in $2k$ variables $t_1,\ldots, t_k, t_{k+1} , \ldots, t_{2k}$ as follows:  
If 
\[
\Phi(a_{i_1},\ldots, a_{i_h}) = \Phi( a_{i_{h+1}} , \ldots, a_{i_{2h}} )
\]
then we have the equation 
\[
\Phi(t_{i_1},\ldots, t_{i_h}) -  \Phi( t_{i_{h+1}} , \ldots, t_{i_{2h}} )= 0.
\]
If 
\[
\Phi(a_{i_1},\ldots, a_{i_h}) < \Phi( a_{i_{h+1}} , \ldots, a_{i_{2h}} )
\]
then 
we have the inequality
\[
\Phi(t_{i_1},\ldots, t_{i_h}) -  \Phi( t_{i_{h+1}} , \ldots, t_{i_{2h}} ) < 0.  
\]
If 
\[
\Phi(a_{i_1},\ldots, a_{i_h}) > \Phi( a_{i_{h+1}} , \ldots, a_{i_{2h}} )
\]
then we  have the inequality
\[
\Phi(t_{i_1},\ldots, t_{i_h}) -  \Phi( t_{i_{h+1}} , \ldots, t_{i_{2h}} ) > 0.  
\]
This procedure gives a system of $k^{2h}$ equations and inequalities 
with rational coefficients for which the set $A  = \{a_1,\ldots, a_k\}$ is a solution in real numbers.  
It follows that the system also has a solution $\hat{B} = \{ \hat{b}_1, \ldots, \hat{b}_k \}$ 
in rational numbers.  

Multiplying an equation or an inequality by a positive number  preserves 
the equation or inequality.  
If $m$ is a positive common multiple of the denominators of the rational numbers 
in the set $\hat{B}$, then $B =  m\ast \hat{B} = \{m \hat{b}_1, \ldots, m\hat{b}_k \}$ 
is a finite set of integers that also solves the system of equations and inequalities.  
Thus, the  set $B$ has the property that, 
for all $h$-tuples $(a_{i_1},\ldots, a_{i_h}) , ( a_{i_{h+1}} , \ldots, a_{i_{2h}} )\in A^h$, we have 
\[
\Phi(a_{i_1},\ldots, a_{i_h}) = \Phi( a_{i_{h+1}} , \ldots, a_{i_{2h}} )
\]
if and only if 
\[
\Phi(b_{i_1},\ldots, b_{i_h}) =  \Phi( b_{i_{h+1}} , \ldots, b_{i_{2h}} )
\]
and so the real set $A$ and the integer set $B$ are $\Phi$-isomorphic.  
Because translation of a set is a $\Phi$-isomorphism, 
we obtain a solution set of positive integers.  
This completes the first proof of Theorem~\ref{MPTQ:theorem:main}.

The second proof uses diophantine approximation.  
We can assume that the linear form $\Phi$ is nonzero with integer coefficients, and so  
\[
\varphi^* = \max(  |\varphi_j| : j=1, \ldots, h) \geq 1.  
\] 
The image under $\Phi$ of a finite set $A$ of real numbers is the set   
\begin{align*}
\Phi(A) 
& = \left\{ \Phi(a_{i_1}, \ldots, a_{i_h}) : (a_{i_1}, \ldots, a_{i_h}) \in A^h  \right\}  \\
& = \left\{ \sum_{j=1}^h  \varphi_j a_{i_j} : a_{i_j} \in A \text{ for } j=1,\ldots, h  \right\}.
\end{align*}
We have  $|\Phi(A)| \geq 2$ because $|A| \geq 2$ and $\Phi$ is nonzero.  
It follows that   
\[
\delta^* = \min \left( | x - x'| : 
x, x' \in \Phi(A) \text{ and } x \neq x'  \right) > 0.  
\]
Let 
\[
0 < \varepsilon < \frac{\min(\delta^*,1)}{2h \varphi^*}.  
\]
By Dirichlet's theorem (Hardy-Wright~\cite[Theorem 201]{hard-wrig08}),  
for every $\varepsilon > 0$, there is a positive integer $q$ 
and a set 
\[
B = \{b_1, \ldots, b_k\}
\]
of integers such that 
\[
|qa_i - b_i| < \varepsilon 
\]
for all $i = 1, \ldots, k$.  We define 
\[
\theta_i = qa_i - b_i 
\]
for $i = 1, \ldots, k$.  
For every $h$-tuple $(a_{i_1},\ldots, a_{i_h}) \in A^h$, we have 
\beq              \label{MPTQ:theta}
\left|   \Phi(\theta_{i_1},\ldots, \theta_{i_h}) \right|  \leq \sum_{j=1}^h |\varphi_j | | \theta_{i_j} | 
 \leq  h \varphi^* \varepsilon   
  <    \frac{ \min \left( \delta^*, 1 \right) }{2}.
\eeq
Define the functions $f: A \rightarrow B$ and $f^*: A^h \rightarrow B^h$ by 
\[
f(a_i) = b_i \qqand f^* (a_{i_1},\ldots, a_{i_h}) = (b_{i_1},\ldots, b_{i_h}).  
\]
Let $(a_{i_1},\ldots, a_{i_h}), (a_{i_{h+1}},\ldots, a_{i_{2h}})$ be $h$-tuples  in $A^h$.  
Using the linearity of $\Phi$, we have  
\[
\Phi(a_{i_1},\ldots, a_{i_h})   = \Phi (a_{i_{h+1}},\ldots, a_{i_{2h}})
\]
if and only if 
\begin{align*}
 \Phi(b_{i_1} , \ldots, b_{i_h})  +  \Phi( \theta_{i_1}, \ldots,  \theta_{i_h}) 
& = \Phi(b_{i_1} + \theta_{i_1}, \ldots, b_{i_h} + \theta_{i_h}) \\
&=  \Phi(qa_{i_1},\ldots, qa_{i_h})  \\
&=  q\Phi(a_{i_1},\ldots, a_{i_h})  \\
&=  q\Phi (a_{i_{h+1}},\ldots, a_{i_{2h}}) \\
& = \Phi (qa_{i_{h+1}},\ldots, qa_{i_{2h}})  \\
& = \Phi (b_{i_{h+1}} + \theta_{i_{h+1}},\ldots, b_{i_{2h}} + \theta_{i_{2h}})  \\
& = \Phi (b_{i_{h+1}},\ldots, b_{i_{2h}}) +  \Phi (\theta_{i_{h+1}},\ldots, \theta_{i_{2h}})
\end{align*}
if and only if 
\[
\Phi(b_{i_1} , \ldots, b_{i_h}) -  \Phi (b_{i_{h+1}},\ldots, b_{i_{2h}}) 
=   \Phi (\theta_{i_{h+1}},\ldots, \theta_{i_{2h}}) -  \Phi( \theta_{i_1}, \ldots,  \theta_{i_h}).
\]
Recall inequality~\eqref{MPTQ:theta}.  
Because $\Phi(b_{i_1} , \ldots, b_{i_h})$ and  $\Phi (b_{i_{h+1}},\ldots, b_{i_{2h}}) $ 
are integers, if $\Phi(b_{i_1} , \ldots, b_{i_h}) \neq  \Phi (b_{i_{h+1}},\ldots, b_{i_{2h}}) $, then 
\begin{align*}
1 & \leq |\Phi(b_{i_1} , \ldots, b_{i_h}) - \Phi (b_{i_{h+1}},\ldots, b_{i_{2h}})  |  \\
&  \leq  \left|  \Phi (\theta_{i_{h+1}},\ldots, \theta_{i_{2h}})  \right| + \left| \Phi( \theta_{i_1}, \ldots,  \theta_{i_h}) \right| \\
& <1
\end{align*}
which is absurd.  
Therefore, $\Phi(a_{i_1} , \ldots, a_{i_h}) = \Phi (a_{i_{h+1}},\ldots, a_{i_{2h}}) $ implies 
$\Phi(b_{i_1} , \ldots, b_{i_h}) =  \Phi   (b_{i_{h+1}},\ldots, b_{i_{2h}}) $, 
and the function $f:A \rightarrow B$ is a $\Phi$-homomorphism.  

Conversely, 
\[
\Phi(b_{i_1} , \ldots, b_{i_h}) = \Phi  (b_{i_{h+1}},\ldots, b_{i_{2h}}) 
\] 
if and only if 
\[
\Phi(qa_{i_1}-\theta_{i_1}, \ldots,qa_{i_h}-\theta_{i_h} ) 
= \Phi(qa_{i_{h+1}}-\theta_{i_{h+1}}, \ldots,qa_{i_{2h}}-\theta_{i_{2h}} ) 
\] 
if and only if 
\[
q (\Phi(a_{i_1} , \ldots, a_{i_h}) - \Phi (a_{i_{h+1}},\ldots, a_{i_{2h}}) )
= \Phi( \theta_{i_1}, \ldots, \theta_{i_h} ) - \Phi (\theta_{i_{h+1}},\ldots, \theta_{i_{2h}}).
\] 
Because $q$ is a positive integer, 
if $ \Phi(a_{i_1} , \ldots, a_{i_h}) \neq \Phi (a_{i_{h+1}},\ldots, a_{i_{2h}}) $, then 
\begin{align*}
\delta^*  & \leq |\Phi(a_{i_1} , \ldots, a_{i_h}) - \Phi (a_{i_{h+1}},\ldots, a_{i_{2h}})  | \\
& \leq q \left| \Phi(a_{i_1} , \ldots, a_{i_h}) - \Phi (a_{i_{h+1}},\ldots, a_{i_{2h}}) \right| \\
&  \leq  \left| \Phi( \theta_{i_1}, \ldots, \theta_{i_h} ) \right| 
+ \left| \Phi  (\theta_{i_{h+1}},\ldots, \theta_{i_{2h}}) \right| \\
& < \min(\delta^*,1) \leq \delta^*
\end{align*}
which is absurd.  
Therefore, $\Phi(b_{i_1} , \ldots, b_{i_h}) = \Phi (b_{i_{h+1}},\ldots, b_{i_{2h}})$ 
implies $\Phi(a_{i_1} , \ldots, a_{i_h}) = \Phi (a_{i_{h+1}},\ldots, a_{i_{2h}})$, 
and the function $f$ is a $\Phi$-isomorphism.

By translation of $B$,  we obtain a set  of positive integers 
that is $\Phi$-isomorphic to $A$.  
This completes the second proof.  

The third proof is the simplest.   
Let $A$ be a nonempty finite set of real numbers that is not a set of integers.    
The additive group $G$ generated by $A$ 
is a torsion-free finitely generated abelian group, and so $G$ is a free abelian group 
of rank $d$ for some positive integer $d$.  Let $f:G \rightarrow \Z^d$ be a group isomorphism.
The restriction of the function $f$  to $A$ is a $\Phi$-isomorphism from $A$ 
to a finite set of lattice points in $\Z^d$.  
By Theorem~\ref{MPTQ:theorem:lattice}, 
this set of lattice points is $\Phi$-isomorphic to a finite set of positive integers.  
This completes the proof.  
\end{proof}

\section{MSTD sets of real numbers and of integers}
In the special case $h=2$ and the linear form 
\[
\Phi(t_1, t_2) = t_1 + t_2
\]
a $\Phi$-isomorphism is called a \textit{Freiman isomorphism}.
The construction of Theorem~\ref{MPTQ:theorem:UpsilonPhi-J} 
applied to the form $\Phi(t_1, t_2) = t_1 + t_2$, with $J = \{ 2\} \subseteq H = \{1,2\}$, 
gives the linear form 
\[
\Phi_J(t_1, t_2) = t_1 - t_2.  
\]
For every subset $A$ of an additive abelian group, we have the 
sumset $A+A = \Phi(A)$ and the difference set $A - A = \Phi_J(A)$.  

\bt                    \label{MPTQ:theorem:main-Z}
Every finite set $A$ of real numbers is Freiman isomorphic to a set $B$ of positive integers. 
If $A$ is an MSTD set of real numbers, then $B$ is an MSTD set of positive integers.  
\et

\begin{proof}
Applying Theorem~\ref{MPTQ:theorem:main} to the linear form $\Phi(t_1, t_2) = t_1 + t_2$, 
we see that, for every finite set $A = \{a_1,\ldots, a_k\}$ of real numbers with $|A| = k$, 
there is a set $B = \{ b_1,\ldots, b_k \}$ of positive integers 
and a function $f:A \rightarrow B$ with  $f(a_i) = b_i$ for $i=1,\ldots, k$,  
such that, for all $a_{i_1}, a_{i_2}, a_{i_3}, a_{i_4} \in A$,
\[
a_{i_1} + a_{i_2} =  a_{i_3} + a_{i_4}
\]
if and only if 
\[
b_{i_1} + b_{i_2} =  b_{i_3} + b_{i_4}.
\]  
Thus, $A$ and $B$ are Freiman isomorphic.  

Using Theorem~\ref{MPTQ:theorem:UpsilonPhi-J} 
with $H = \{1,2\}$ and  $J = \{2\}$, we have 
\[
a_{i_1} - a_{i_2} =  a_{i_3} - a_{i_4}
\]
if and only if 
\[
b_{i_1} - b_{i_2} =  b_{i_3} - b_{i_4}.
\]
It follows that 
\[
|A + A|  = |B + B| 
\]
and
\[
|A - A|  = |B - B|
\]
and so $B$ is an MSTD set of integers if $A$ is an MSTD set of real numbers.  
\end{proof}

\bt                \label{MPTQ:theorem:main-real8}
If $A$ is an MSTD set of real numbers, then $|A| \geq 8$.
If $A$ is an MSTD set of real numbers and $|A| = 8$, 
then $A$ is Freiman isomorphic to the set $\{0, 2, 3, 4, 7, 11, 12, 14\}$.
\et

\begin{proof}
This follows immediately from Hegarty's theorem that, 
if $B$ is an MSTD set of integers, then $|A| \geq 8$, 
and if $B$ is an MSTD set of integers and $|A| = 8$, 
then $B$ is affinely isomorphic to the set $\{0, 2, 3, 4, 7, 11, 12, 14\}$.
\end{proof}

An additive abelian group $G$ is \emph{2-divisible} if for every $a \in G$ 
there exists $x\in G$ such that $2x=a$.  We write $x = a/2$.  
For example, the group $\Z/3\Z$ is 2-divisble.

\bt                            \label{MPTQ:theorem:affineMap}
Let $A^* = \{0, 2, 3, 4, 7, 11, 12, 14\}$, and let $B $ be a subset of a 2-divisible group 
such that $B$ is Freiman isomorphic to $A$.  
If $f: A^* \rightarrow B$ is a  Freiman isomorphism, then $f$ is the affine map 
\beq                      \label{MPTQ:affineMap}
f(x) =  x \left(\frac{f(2)-f(0)}{2} \right) + f(0) 
\eeq
for all $x \in A^*$. 
\et

\begin{proof}
Because $f$ is a Freiman isomorphism, 
\[
4+ 0 = 2 + 2
\]
implies that 
\[
 f(4) + f(0)= f(2) + f(2) 
\]
and so
\[
f(4) = 2f(2) - f(0) = \frac{4f(2) - 2f(0)}{2}.
\]
Similarly, $3+3 = 2+4$ implies that 
\[
f(3) = \frac{f(2) + f(4)}{2} = \frac{3f(2) - f(0)}{2}.
\]
The equation $7+0 = 4+3$ implies that 
\[
f(7) =  \frac{7f(2) - 5f(0)}{2}.
\]
The equation $11+0 = 7+4$ implies that 
\[
f(11) =  \frac{11f(2) - 9f(0)}{2}.
\]
The equation $14+0 = 7 + 7$ implies that 
\[
f(14) =  \frac{14f(2) - 12f(0)}{2}.
\]
The equation  $12 + 2 = 14+0$  implies that 
\[
f(12) =  \frac{12f(2) - 10f(0)}{2}.
\]
We see that  
\[
f(x) = \frac{xf(2) - (x-2)f(0)}{2} = x \left(\frac{f(2)-f(0)}{2} \right)  + f(0) 
\]
for all $x \in A^*$. 
This completes the proof.  
\end{proof}

The following result is due to Moshe Newman~\cite{newm16}.

\bt[Newman]                                \label{MPTQ:theorem:MosheNewman}
Every MSTD set $A$ of real numbers with $|A| = 8$  
is affinely isomorphic to the set $A^* = \{0, 2, 3, 4, 7, 11, 12, 14\}$.
\et

\begin{proof}
By Theorem~\ref{MPTQ:theorem:main-real8}, 
the set $A$ of  real numbers  is Freiman isomorphic to $A^*$.  
Let $A = \{ a_1, a_2, \ldots, a_8\}$, where $a_1 < a_2 < \cdots < a_8$. 
Define the affine map $g:\R \rightarrow \R$ by 
\[
g(x) = \frac{2(x-a_1)}{a_2 - a_1}.
\]
Let $b_i = g(a_i)$ for $i = 1, \ldots, 8$, and let $B = g(A) = \{b_1, \ldots, b_8\}$.  
Note that $b_1 = f(a_1) = 0 $ and $b_2 = f(a_2) = 2$, 
and that $b_1 < b_2 < \cdots < b_8$.  
The set  $B$ of real numbers is affinely isomorphic to $A$ 
and Freiman isomorphic to  $A^*$.  
Let 
\[
f:A^* \rightarrow B 
\]
be a Freiman isomorphism.  
By Theorem~\ref{MPTQ:theorem:affineMap}, the function $f$ is an affine map 
of the form~\eqref{MPTQ:affineMap}.  
If $f(2) > f(0)$, then $f$ is strictly increasing, and so $f(0) = 0$ and $f(2) = 2$.   
This implies that $f(x) = x$ for all $x \in A^*$, and so $A^* = B$ is affinely isomorphic to $A$.

If $f(2) < f(0)$, then $f$ is strictly decreasing, and $f(0) = 14$ and $f(2) = 12$.
This implies that $f(x) = 14-x$ for $x \in A^*$, and so $B = 14 - A^* = \{0,2,3,7, 10,11,12,14\}$,  
which is also affinely isomorphic to $A^*$.
This completes the proof.  
\end{proof}

\section{More products than quotients}
Let $B$ be a set of positive real numbers.  We define the \textit{product set}
\[
B \cdot B = \{ b \cdot b':b,b' \in B\}
\]
and the quotient set
\[
B/B = \{ b/b':b,b' \in B \}.  
\]

The observation that 
\[
3 \cdot 2 = 2 \cdot 3 
\]
but
\[
\frac{3}{2} \neq \frac {2}{3}
\]
suggests that a finite set of positive real numbers 
should have more quotients, or at least as many quotients, as products, 
but it easy to show that there do exist sets with more products than quotients.  
Such sets are called \textit{MPTQ sets}.  

Here is a simple construction of MPTQ sets.
Let $A$ be a set of real numbers.  For every positive real number $c$,  let 
\[
c^A = \{c^a:a \in A\}.
\]
The function $a \mapsto c^a$ is a Freiman isomorphism from 
the additive set $A$ to the multiplicative set $c^A$.
The set  $A$ is an MSTD set of real numbers  
if and only if  $c^A = \{c^a:a \in A\}$ is an MPTQ set of positive real numbers.  
Moreover, if $A$ is an MSTD set of positive integers and if $c > 1$ is an integer, 
then $c^A$ is an MPTQ set of positive integers.  

Conversely, let $B$ be a set of positive real numbers.  
For every positive real number $c$, let 
\[
\log_c B = \left\{ \log_c b: b \in B \right\}.
\]
The function $b \mapsto  \log_c b$ is a Freiman isomorphism from 
 the multiplicative set $B$ to the additive set $\log_c B$.
The set  $B$ is an MPTQ set of positive real numbers 
if and only if $\log_c B$  is an MSTD set of real numbers.

\section{Problems}

Hegarty gave a computational proof of the theorem that an MSTD set of integers 
must contain at least 8 elements.   
A large number of cases must be checked, and a computer checks them.  
This calculation does not explain why the result is true.  
\bprob
Is there a proof of Hegarty's theorem that explains why an MSTD set of integers 
cannot have 7 elements?  
\eprob

A subset $A$ of an abelian group $G$ is \textit{symmetric}
 if there exists an element $s \in G$ such 
that $A =\{s-a:a\in A\}$.  
For example, every finite arithmetic progression is symmetric.

\bl              \label{MPTQ:lemma:symmetricForm}
Let 
\[
\Phi(t_1,\ldots, t_h) = \sum_{ j =1}^h \varphi_j t_j
\]
be a nonzero linear form with integer coefficients, 
and, for $J \subseteq H = \{1,2,\ldots, h\}$, let 
\[
\Phi_J(t_1,\ldots, t_h) = \sum_{j \in H\setminus J} \varphi_j t_j - \sum_{j \in J} \varphi_j t_j.  
\]
If $A$ is a symmetric finite subset of an abelian group $G$, then 
\[
|\Phi(A)| = |\Phi_J(A)|.
\] 
\el

\begin{proof}
Because $A$ is symmetric, there is an element $s \in G$ such that, for every $a \in A$ 
there is a unique $a' \in A$ with $a = s-a'$.  Let 
\[
s^* = \sum_{j \in J} \varphi_j  s 
\]
If $x \in \Phi(A)$, then there exists $(a_1,\ldots, a_h) \in A^h$ such that 
\begin{align*}
x & = \Phi(a_1,\ldots, a_h) \\
& =  \sum_{j \in H\setminus J} \varphi_j a_j + \sum_{j \in J} \varphi_j a_j \\
& =  \sum_{j \in H\setminus J} \varphi_j a_j + \sum_{j \in J} \varphi_j (s-a'_j) \\
& =  \sum_{j \in J} \varphi_j  s
+  \sum_{j \in H\setminus J} \varphi_j a_j - \sum_{j \in J} \varphi_j a'_j \\
& \in s^* + \Phi_J(A)
\end{align*}
and so $\Phi(A) \subseteq \{s^*\} + \Phi_J(A)$  and $|\Phi(A)| \leq |\Phi_J(A)|$.  
Similarly, $\Phi_J(A) \subseteq \{-s^*\} + \Phi(A)$ and $|\Phi_J(A)| \leq |\Phi(A)|$.  
This completes the proof.  
\end{proof}

Here is a simple problem chosen from a large class of related problems 
about linear forms and finite sets of integers.  
Consider the linear forms 
\[
\Phi(t_1,t_2,t_3) = t_1 + t_2 + t_ 3
\]
and
\[
\Phi_3(t_1,t_2,t_3) = t_1 + t_2 - t_ 3
\]
For every set $A$ of integers, we have 
\[
\Phi(A) = A+A+A \qqand \Phi_3(A) = A+A-A.
\]
By Lemma~\ref{MPTQ:lemma:symmetricForm}, 
if $A$ is a finite symmetric set of integers, then $|\Phi(A)| = |\Phi_3(A)| $.
For example, if $A$ is an arithmetic progression of length $k$, then 
$|\Phi(A)| = |\Phi_3(A)| = 3k-2$.  
If $A =  \{0, 2, 3, 4, 7, 11, 12, 14\}$, then 
\[
\Phi(A) = \Phi_3(A) + \{14\} = [0,42] \setminus \{1,41\} 
\]
and 
\[
|\Phi(A)| = |\Phi_3(A)| = 41.  
\]
The symmetry of the form $\Phi$ and the asymmetry of the form $\Phi_3$ 
suggest that $|\Phi(A)| \leq |\Phi_3(A)|$ for most finite sets of integers.  
\bprob
Describe the finite sets $A$ of integers such that $|\Phi(A)| > |\Phi_3(A)| $.
\eprob

\bprob
Can we go beyond linear forms?   
Consider quadratic, cubic, and other higher degree forms, 
or general polynomials in $h$ variables, such as the polynomials $t_1 + t_2^2$
and $t_1^2 + t_2^2$.  
Here is a sample problem.  To every $h$-tuple $I = (i_1,\ldots, i_h) $ 
of nonnegative integers, we associate the monomial $t^I = t_1^{i_1} \cdots t_h^{i_h}$.  
Let $\mci$ be a finite set of $h$-tuples, and let 
\[
P(t_1,\ldots, t_h) = \sum_{I \in \mci} \varphi_I t^I 
\]
be a polynomial with nonzero integer coefficients.  Let $\varepsilon_I \in \{\pm 1\}$ for $I \in \mci$, 
and let $\varepsilon = (\varepsilon_I)_{I \in \mci}$.  
Define 
\[
P_{\varepsilon}(t_1,\ldots, t_h) = \sum_{I \in \mci}\varepsilon_I  \varphi_I t^I.
\]
Compare the cardinalities of the sets $P_{\varepsilon}(A)$ for finite sets $A$ of integers.
\eprob

\bprob
How can one model the behavior of polynomial images of finite sets 
of real numbers by finite sets of integers?  
 \eprob

\section{Acknowledgements}
This paper originated in  a conference in honor of Grigori Freiman 
at Tel Aviv University in July, 2016. 
I stated in my lecture that I did not know who had discovered the first example of
an MSTD set, but that it might have been a  set of real numbers.
Noga Alon remarked that, by the method in his paper~\cite{alon-klei90} with Kleitman,
one can always construct an MSTD set of integers from an MSTD set of
real numbers.  This gave the first proof 
of Theorem~\ref{MPTQ:theorem:main}. 
 Kevin O'Bryant suggested the use of Dirichlet's theorem in the diophantine approximation proof 
 of Theorem~\ref{MPTQ:theorem:main}.  I described both proofs in a lecture at the Integers Conference in October, 2016, and Paul Pollack 
 observed that the structure theorem for finitely generated torsion-free abelian groups 
 gives the third proof.  
 
 Theorems~\ref{MPTQ:theorem:affineMap} and~\ref{MPTQ:theorem:MosheNewman} 
 are essentially due to Moshe Newman~\cite{newm16}.

\def\cprime{$'$} \def\cprime{$'$} \def\cprime{$'$} \def\cprime{$'$}
\providecommand{\bysame}{\leavevmode\hbox to3em{\hrulefill}\thinspace}
\providecommand{\MR}{\relax\ifhmode\unskip\space\fi MR }
\providecommand{\MRhref}[2]{%
  \href{http://www.ams.org/mathscinet-getitem?mr=#1}{#2}
}
\providecommand{\href}[2]{#2}

\end{document}